# VISCOSITY SOLUTIONS OF HAMILTON–JACOBI EQUATIONS WITH DISCONTINUOUS COEFFICIENTS

GIUSEPPE MARIA COCLITE AND NILS HENRIK RISEBRO

Abstract. We consider Hamilton–Jacobi equations, where the Hamiltonian depends discontinuously on both the spatial and temporal location. Our main results are the existence and well–posedness of a viscosity solution to the Cauchy problem. We define a viscosity solution by treating the discontinuities in the coefficients analogously to "internal boundaries". By defining an appropriate penalization function, we prove that viscosity solutions are unique. The existence of viscosity solutions is established by showing that a sequence of front tracking approximations is compact in $L^\infty$, and that the limits are viscosity solutions.

## 1. Introduction

In this paper we study the initial value problem

$$\begin{cases} u_t + H\left(u_x, a(x), g(t)\right) = 0 & x \in \mathbb{R}, \quad t > 0, \\ u(x,0) = u_0(x) & t = 0. \end{cases} \tag{1.1}$$

where $H : \mathbb{R}^3 \to \mathbb{R}$. The functions $a : \mathbb{R} \to \mathbb{R}$ and $g : \mathbb{R}^+ \to \mathbb{R}$ are called "coefficients", and are allowed to be discontinuous. The classical theory for viscosity solutions of Hamilton–Jacobi equations, see e.g. [4] does not cover the case where the Hamiltonian $H$ is discontinuous. This is because the straightforward method of comparing sub- and supersolutions does not work if $H$ is discontinuous in $x$ or $t$.

One application of Hamilton–Jacobi equations with discontinuous coefficients is the optical shape-from-shading problem

$$\frac{1}{\sqrt{1+u_x^2+u_y^2}} = I(x,y), \tag{1.2}$$

where $I$ denotes the intensity of the reflected light, and $u = u(x,y)$ the height of the underlying surface. If the gradient of $u$ is discontinuous, then the intensity of the reflected light will vary discontinuously. Another related application is the synthetic radar shape-from-shading equation

$$\frac{u_x^2}{\sqrt{1+u_x^2+u_y^2}} = I(x,y). \tag{1.3}$$

Also in this case, if the gradient of $u$ is discontinuous, so is $I$. Note that (1.3) can be reformulated as the evolution equation

$$u_x - I(x,y)\sqrt{\frac{1}{2} + \sqrt{\frac{1}{4} + \frac{1+u_y^2}{I^2(x,y)}}} = 0.$$

This is reminiscent of (1.1), but unless $I(x,y) = f(a(y), g(x))$ for some function $f$, this equation is not actually of the type considered in this paper.

A number of authors have considered Hamilton–Jacobi equations with coefficients that have some kind of singularity. In [9] Ishii and Ramaswami considered a boundary, and an initial value problem. In this paper [9] the notion of viscosity solution was extended to some classes of

*Date*: September 17, 2018.

The research was funded [in part] by the BeMatA program of the Research Council of Norway and by the European network HYKE, funded by the EC as contract HPRN-CT-2002-00282.

discontinuous Hamiltonians. Using the definitions from [9], in Tourin [21] showed that for this class of discontinuous Hamiltonians, a unique viscosity solution exists.

Also, if the Hamiltonian $H$ satisfies some structure conditions, essentially amounting to discontinuous jumps "in one direction", there exists a unique solution to the initial value problem, see Capuzzo Dolcetta and Perthame [2] or Strömberg [19].

For shape-from-shading problems, Ostrov [16, 17] showed that provided the viscosity solution was bounded and absolutely continuous, then a sequence of approximations, defined by smoothing the coefficients, converged to the unique viscosity solution. This viscosity solution was defined as the solution of an auxiliary control problem. The technique of using the associated control problem (for convex Hamiltonians) was also exploited by Dal Maso and Frankowska in [5, 6] when studying some Hamilton–Jacobi equations with a discontinuous Hamiltonian.

The Hamilton–Jacobi equation (1.1) is formally equivalent to the conservation law

$$p_t + H(p, a, g)_x = 0, \tag{1.4}$$

where $p = u_x$. This is an example of a conservation law with discontinuous coefficients, and such equations have been extensively studied by a number of authors, see e.g. Klingenberg and Risebro [13, 14], Gimse and Risebro [8], Klaussen and Risebro [12]. In [12] it was shown that smoothing the coefficients for the conservation law produced a convergent sequence of solutions, such that the limit was a weak solution of the conservation law. For a class of flux functions, Towers [22] proved uniqueness within the class of piecewise smooth solutions by using a variant of the Kružkov [15] approach, as well as convergence of monotone difference methods [23]. Finally, Seguin and Vovelle [18] studied a special case of the purely hyperbolic version of (1.4) with the flux fuction taking the form $H(a, p) = ag(p)$ The authors proved uniqueness of $L^\infty$ entropy solutions by the Kružkov method [15]. Recently, a quite general theory for conservation laws with discontinuous fluxes was established by Karlsen, Risebro and Towers in [11], see also Coclite and Risebro [3] in which well-posedness is established for some conservation laws of the type (1.4).

In this paper we establish both uniqueness and existence of viscosity solutions to (1.1). The uniqueness is demonstrated by using a standard "doubling of the variables" approach, and choosing the penalization function so that the discontinuities does not influence the technique. This is motivated by the techniques used in [11], although we could possibly avoided some of the technicalities by using a related approach used by Benth, Karlsen and Reikvam in [1], for a control problem.

The existence is established by showing that a front tracking scheme, designed for (1.4) is well-defined and produces a sequence that converges to an appropriate entropy solution. The integrals of the approximate solutions to the conservation law, are then shown to be approximate solutions to the Hamilton–Jacobi equation, and their limit is a viscosity solution. This program was carried out for conservation laws/Hamilton–Jacobi equations without $x$ or $t$ dependence by Karlsen and Risebro in [10].

To be specific, we make the assumptions:

(**A.1**) The Hamiltonian $H$ satisfies the following Lipschitz conditions

$$\left|H(p, a, g) - H(p, a', g')\right| \le C\left(|g - g'| + |a - a'|\right)(1 + |p|), \tag{1.5}$$

$$\left|H(p, a, g) - H(q, a, g)\right| \le C\left|p - q\right|, \tag{1.6}$$

for all $p$, $q$, $a$, $a'$, $g$, $g'$, and some positive constant $C$.

(**A.2**)

$$\frac{\partial H}{\partial a} \ge 0, \quad \text{and} \quad \frac{\partial H}{\partial g} \ge 0. \tag{1.7}$$

(**A.3**)

$$\frac{\partial H}{\partial p}(0, a, g) = 0 \quad \text{and} \quad \frac{\partial^2 H}{\partial p^2}(0, a, g) < 0, \tag{1.8}$$

for all $a$ and $g$.

(**A**.4) The map $p \mapsto H(p,a,g)$ is even, and strictly increasing in the interval $\langle -\infty, 0\rangle$, and strictly decreasing in the interval $\langle 0,\infty\rangle$, for all $a$ and $g$.

(**A**.5) We have that

$$\lim_{p\to\infty} \frac{|H(p,a,g)|}{p} \geq C > 0, \tag{1.9}$$

for some constant $C$ that is independent of $a$ and $g$.

(**A**.6) The coefficient $a(x)$ is piecewise continuously differentiable, with finitely many jump discontinuities in $a$ and $a'$, located at the points $x_1 < \cdots < x_M$.

(**A**.7) The coefficient $g(t)$ is piecewise continuous, with finitely many jump discontinuities in $g$ and $g'$, located at the points $t_1 < t_2 < \cdots < t_N$.

(**A**.8) The coefficients $a$ and $g$ are of bounded variation, i.e.,

$$|g|_{BV(\mathbb{R}^+)} < \infty, \quad |a|_{BV(\mathbb{R})} < \infty \tag{1.10}$$

The coefficients have bounded derivatives away from discontinuities, i.e.,

$$\max_{x\notin\{x_1,\dots,x_M\}} |a'(x)| < \infty, \qquad \max_{t\notin\{t_1,\dots,t_N\}} |g'(t)| < \infty. \tag{1.11}$$

A viscosity solution of the Cauchy problem (1.1) is defined as follows:

**Definition 1.1.** *Fix $T > 0$, and let $u : \Pi_T = \mathbb{R}\times[0,T\rangle \to \mathbb{R}$ be a bounded and uniformly continuous (BUC) function. We call $u$ a viscosity subsolution of* (1.1) *if $u(x,\cdot) = u_0$, and if for all $C^1$ functions $\varphi(x,t)$ such that $u-\varphi$ has a local maximum at $(x_0,t_0)\in\Pi_T$, then*

$$\varphi_t(x_0,t_0) + \min\left\{H\left(\varphi_x(x_0,t_0), a\left(x_0^-\right), g\left(t_0^-\right)\right), H\left(\varphi_x(x_0,t_0), a\left(x_0^+\right), g\left(t_0^-\right)\right)\right\} \leq 0. \tag{1.12}$$

*Analogously, we call $u$ a viscosity supersolution if $u(x,\cdot) = u_0$, and if for every $C^1$ function $\varphi$ such that $u-\varphi$ has a local minimum at $(x_0,t_0)$, then*

$$\varphi_t(x_0,t_0) + \max\left\{H\left(\varphi_x(x_0,t_0), a\left(x_0^-\right), g\left(t_0^-\right)\right), H\left(\varphi_x(x_0,t_0), a\left(x_0^+\right), g\left(t_0^-\right)\right)\right\} \geq 0. \tag{1.13}$$

*Finally, we say that $u$ is a viscosity solution of* (1.1) *if it is both a viscosity subsolution and a supersolution.*

Setting $p = u_x$, it formally follows that $p$ solves the Cauchy problem

$$\begin{cases} p_t + H(p,a(x),g(t))_x = 0 & x\in\mathbb{R}, \quad t>0,\\ p(x,0) = p_0(x) & t = 0, \end{cases} \tag{1.14}$$

where $H$, $a$ and $g$ are as before. To define an entropy weak solution of (1.14) we shall need the Temple singular mapping $\Psi$ defined by

$$\Psi(p,a,g) = \operatorname{sign}(p)\,\frac{H(p,a,g) - H(0,a,g)}{H(0,a,g)}. \tag{1.15}$$

By an entropy solution to (1.14) we shall mean a function $p$ satifying the following definition:

**Definition 1.2.** *Let $p : \Pi_T \to \mathbb{R}$ be a measurable function. We say that $p$ is an entropy solution of* (1.14) *if the following hold:*

(**D**.1) *$p\in L^1(\Pi_T)\cap L^\infty(\Pi_T)$, $\Psi(p(\cdot,t),a,g(t))\in BV(\mathbb{R})$ for all $t\in[0,T\rangle$. Furthermore, the map $[0,T\rangle \ni t\mapsto p(t,\cdot)$ is $L^1(\mathbb{R})$ Lipschitz continuous.*

(**D**.2) *The function $p$ is a weak solution of* (1.14), *i.e.,*

$$\iint_{\Pi_T} p\varphi_t + H(p,a,g)\varphi_x \,dtdx + \int_{\mathbb{R}} \varphi(x,0)p_0\,dx = 0, \tag{1.16}$$

*for all test functions $\varphi\in C_0^1(\Pi_T)$.*

(**D.3**) *The following (entropy) inequality holds for all constants $c$ and all non-negative test functions $\varphi$,*

$$\iint_{\Pi_T} |p-c|\,\varphi_t + F(p,x,t,c)\varphi_x \,dtdx - \sum_{m=0}^{M} \int_{x_m}^{x_{m+1}} \int_0^T \operatorname{sign}(p-c)\, H_a(c, a(x), g(t)) a'(x) \varphi \,dtdx$$
$$+ \sum_{m=1}^{M} \int_0^T \left| H(c, a\left(x_m^+\right), g(t)) - H(c, a\left(x_m^-\right), g(t)) \right| dt \geq 0, \tag{1.17}$$

*where we have set $x_0 = -\infty$, $x_{M+1} = \infty$, and $F$ is given by*

$$F(p,x,t,c) = \operatorname{sign}(p-c)\left[H(p, a(x), g(t)) - H(p, a(x), g(t))\right], \quad t > 0, \quad x \in \mathbb{R}.$$

The main result of this paper is summarized by

**Main Theorem.** *Assume that the assumptions (**A.1**) – (**A.8**) hold. Let $u_0$ and $v_0$ be two functions in $BUC(\mathbb{R})$, then there exist corresponding viscosity solutions of* (1.1) *$u = u(x,t)$ and $v = v(x,t)$, satisfying the initial conditions*

$$u(x,0) = u_0(x), \qquad v(x,0) = v_0(x),$$

*furthermore*

$$\|u(\cdot,t) - v(\cdot,t)\|_{L^\infty(\mathbb{R})} \leq \|u_0 - v_0\|_{L^\infty(\mathbb{R})},$$

*for all $t \in [0,T]$. Moreover the functions $p$ and $q$ defined by*

$$p = \frac{\partial u}{\partial x}, \quad \text{and} \quad q = \frac{\partial v}{\partial x},$$

*are entropy weak solutions of* (1.14)*, taking the initial values*

$$p(\cdot,0) = \frac{\partial u_0}{\partial x}, \quad \text{and} \quad q(\cdot,0) = \frac{\partial v_0}{\partial x}.$$

**Remark 1.3.** *From* [3] *and* [11] *we know that*

$$\|p(\cdot,t) - q(\cdot,t)\|_{L^1(\mathbb{R})} \leq \|p(\cdot,0) - q(\cdot,0)\|_{L^1(\mathbb{R})}.$$

The remainder of this paper is organized as follows. In Section 2 we show that viscosity solutions are unique by establishing a comparison principle. Then in Section 3 we show existence of a solution by first showing that a front tracking algorithm for (1.4) is well defined. Then we show that this automatically yields a front tracking algorithm for the Hamilton–Jacobi equation. The convergence of sequences produced by this algorithm(s) finally gives the existence.

## 2. Uniqueness and stability

This section is devoted to the proof of the comparison principle.

**Theorem 2.1** (**Comparison principle**)**.** *Assume that (**A.1**), (**A.6**), (**A.7**) and (**A.8**) hold. Let $u$ and $v$ be two viscosity solutions of* (1.1) *taking the initial values*

$$u(x,0) = u_0(x), \qquad v(x,0) = v_0(x),$$

*where $u_0$ and $v_0$ are in $BUC(\mathbb{R})$. If*

$$u_0(x) \leq v_0(x), \quad x \in \mathbb{R},$$

*then*

$$u(x,t) \leq v(x,t), \quad \text{for } x \in \mathbb{R} \text{ and } t \in [0,T].$$

An immediate consequence of the comparison principle is that viscosity solutions are unique. Without loss of generality we shall assume that

$$x_M = 0, \tag{2.1}$$

which in particular implies that $a$ is continuous in the interval $\langle 0, \infty\rangle$. To prove Theorem 2.1 we shall need two lemmas, the first of which is

**Lemma 2.2.** *Assume that (**A**.1), (**A**.6), (**A**.7) and (**A**.8) hold. Let $u$ be a viscosity solution of* (1.1) *taking the initial data $u_0 \in BUC(\mathbb{R})$. Let $\varphi$ be a continuously differentiable test function $\varphi = \varphi(x,t)$. If $u - \varphi$ attains a local maximum at a point $(x_0, T)$, where $x_0 \notin \{x_1, \dots, x_M\}$, then*

$$\varphi_t(x_0, T) + H\left(\varphi_x(x_0, T), a(x_0), g\left(T^-\right)\right) \le 0, \tag{2.2}$$

*while if $u - \varphi$ attains a local minimum at $(x_0, T)$ then*

$$\varphi_t(x_0, T) + H\left(\varphi_x(x_0, T), a(x_0), g\left(T^-\right)\right) \ge 0. \tag{2.3}$$

*Proof.* Assume that $u - \varphi$ has a local maximum at $(x_0, T)$. Without loss of generality, we can assume that this maximum is strict. Now for $\varepsilon > 0$, define

$$\varphi_\varepsilon(x,t) = \varphi(x,t) + \frac{\varepsilon}{T-t}.$$

For each $\varepsilon$, $u - \varphi_\varepsilon$ attains a local maximum at a point $(x_\varepsilon, t_\varepsilon)$, such that

$$x_\varepsilon \to x_0, \quad \text{and} \quad t_\varepsilon \to T, \quad \text{as } \varepsilon \to 0.$$

Since $x_0 \notin \{x_1, \dots, x_M\}$ by (**A**.6), for sufficiently small $\varepsilon$, also $x_\varepsilon \notin \{x_1, \dots, x_M\}$, and $t_\varepsilon \notin \{t_1, \dots, t_N\}$. Since $u$ is a viscosity solution,

$$\varphi_{\varepsilon,t}(x_\varepsilon, t_\varepsilon) + H\left(\varphi_{\varepsilon,x}(x_\varepsilon, t_\varepsilon), a(x_\varepsilon), g(t_\varepsilon)\right) \le 0.$$

We also have that

$$\varphi_{\varepsilon,x}(x,t) = \varphi_x(x,t) \quad \text{and} \quad \varphi_{\varepsilon,t}(x,t) = \varphi_t(x,t) + \frac{\varepsilon}{(T-t)^2},$$

and thus

$$\varphi_t(x_\varepsilon, t_\varepsilon) + H\left(\varphi_x(x_\varepsilon, t_\varepsilon), a(x_\varepsilon), g(t_\varepsilon)\right) \le -\frac{\varepsilon}{(T-t_\varepsilon)^2} < 0.$$

Therefore, by (**A**.1), (**A**.6) and (**A**.7), (2.2) holds. The inequality (2.3) is proved by the same arguments. □

The second lemma we shall need to show Theorem 2.1 is

**Lemma 2.3.** *Assume that (**A**.1), (**A**.6), (**A**.7) and (**A**.8) hold. Let $u$ and $v$ be two viscosity solutions of* (1.1) *in the strip $\Pi_S$, where $0 < S < t_1$, taking the initial values*

$$u(x,0) = u_0(x), \qquad v(x,0) = v_0(x),$$

*where $u_0$ and $v_0$ are in $BUC(\mathbb{R})$. If*

$$u_0(x) \le v_0(x), \quad x \in \mathbb{R}, \tag{2.4}$$

*then*

$$u(x,t) \le v(x,t), \quad \textit{for } x \in \mathbb{R} \textit{ and } t \in [0,S].$$

*Proof.* We assume that (2.4) does not hold. Then there exists a constant $\lambda > 0$ such that

$$\sup_{0 \le t \le S,\, x \in \mathbb{R}} \left[u(x,t) - v(x,t) - 2\lambda t\right] \ge 0. \tag{2.5}$$

For $\varepsilon > 0$ define the function

$$\begin{aligned}
\Phi_\varepsilon(x,t,y,s) &= u(x,t) - v(y,s) - \lambda(t+s) - \varepsilon(x^2+y^2) \\
&\quad - \frac{(t-s)^2 + (|x|-|y|)^2}{\varepsilon^2} - \frac{1}{\varepsilon}\left(\frac{1}{x^2} + \frac{1}{y^2}\right) - \frac{\varepsilon^3}{(|x|-|y|)^2 + \varepsilon^2} - K\chi(x,y),
\end{aligned} \tag{2.6}$$

where

$$K = 2\left(\max|u| + \max|v|\right), \quad \text{and} \quad \chi(x,y) = \begin{cases} 1 & \text{if } x \le 0 \text{ or } y \le 0, \\ 0 & \text{otherwise.} \end{cases}$$

The function $\Phi_\varepsilon$ has a global maximum at a point $(x_\varepsilon, t_\varepsilon, y_\varepsilon, s_\varepsilon) \in (\mathbb{R} \times [0,S])^2$, and due to the presence of the term $K\chi$,

$$x_\varepsilon > 0 \quad \text{and} \quad y_\varepsilon > 0.$$

Next we show that we can find two constants $0 < C_1 < C_2$ such that

$$\frac{C_1}{\sqrt{\varepsilon}} \leq x_\varepsilon,\ y_\varepsilon \leq \frac{C_2}{\sqrt{\varepsilon}}, \quad C_1\varepsilon^{3/2} - \varepsilon \leq |x_\varepsilon - y_\varepsilon| \leq C_2\varepsilon, \quad \text{and} \quad |t_\varepsilon - s_\varepsilon| \leq C_2\varepsilon. \tag{2.7}$$

Since

$$2\Phi_\varepsilon(x_\varepsilon, t_\varepsilon, y_\varepsilon, s_\varepsilon) \geq \Phi_\varepsilon(x_\varepsilon, t_\varepsilon, x_\varepsilon, t_\varepsilon) + \Phi_\varepsilon(y_\varepsilon, s_\varepsilon, y_\varepsilon, s_\varepsilon),$$

we find that

$$2\frac{(t_\varepsilon - s_\varepsilon)^2 + (x_\varepsilon - y_\varepsilon)^2}{\varepsilon^2} + 2\frac{\varepsilon^3}{(x_\varepsilon - y_\varepsilon)^2 + \varepsilon^2} \leq u(x_\varepsilon, t_\varepsilon) - v(y_\varepsilon, s_\varepsilon) - u(y_\varepsilon, s_\varepsilon) + v(x_\varepsilon, t_\varepsilon).$$

Since $u$ and $v$ are globally bounded, the terms on the left hand side of the above inequality is also bounded, and thus we can find constants $0 < C_3 < C_4$ such that

$$C_3\varepsilon^{3/2} - \varepsilon \leq |x_\varepsilon - y_\varepsilon| \leq C_4\varepsilon, \quad \text{and} \quad |t_\varepsilon - s_\varepsilon| \leq C_4\varepsilon.$$

Next we observe that

$$\begin{aligned}\Phi_\varepsilon(x_\varepsilon, t_\varepsilon, y_\varepsilon, s_\varepsilon) &\geq \Phi_\varepsilon\left(\frac{1}{\sqrt{\varepsilon}}, t_\varepsilon, \frac{1}{\sqrt{\varepsilon}}, s_\varepsilon\right) \\ &= u\left(\frac{1}{\sqrt{\varepsilon}}, t_\varepsilon\right) - v\left(\frac{1}{\sqrt{\varepsilon}}, s_\varepsilon\right) - \lambda(t_\varepsilon + s_\varepsilon) - 4 - \frac{(t_\varepsilon - s_\varepsilon)^2}{\varepsilon^2} - \varepsilon.\end{aligned}$$

This implies that

$$\begin{aligned}\varepsilon\left(x_\varepsilon^2 + y_\varepsilon^2\right) + \frac{1}{\varepsilon}\left(\frac{1}{x_\varepsilon^2} + \frac{1}{y_\varepsilon^2}\right) \leq u(x_\varepsilon, t_\varepsilon) - v(y_\varepsilon, s_\varepsilon) - u\left(\frac{1}{\sqrt{\varepsilon}}, t_\varepsilon\right) + v\left(\frac{1}{\sqrt{\varepsilon}}, s_\varepsilon\right) \\ - \frac{(x_\varepsilon - y_\varepsilon)^2}{\varepsilon^2} - \frac{\varepsilon^3}{(x_\varepsilon - y_\varepsilon)^2 + \varepsilon^2} + 4 - \varepsilon.\end{aligned}$$

The right hand side of this inequality is bounded independently of $\varepsilon$, and hence we can find constants $0 < C_5 < C_6$ such that

$$\frac{C_5}{\sqrt{\varepsilon}} \leq x_\varepsilon,\ y_\varepsilon \leq \frac{C_6}{\sqrt{\varepsilon}}.$$

This finishes the proof of (2.7). We can extract some more information from this, since

$$\frac{1}{\varepsilon}\left(\frac{1}{x_\varepsilon^2} - \frac{1}{y_\varepsilon^2}\right) = \left(y_\varepsilon^2 - x_\varepsilon^2\right) \cdot \frac{1}{\varepsilon x_\varepsilon^2} \cdot \frac{1}{\varepsilon y_\varepsilon^2} \cdot \varepsilon.$$

Therefore

$$\lim_{\varepsilon \to 0} \frac{1}{\varepsilon}\left(\frac{1}{x_\varepsilon^2} - \frac{1}{y_\varepsilon^2}\right) = 0. \tag{2.8}$$

Next, we aim to show the three following limits

$$\begin{aligned}\lim_{\varepsilon \to 0} \frac{t_\varepsilon - s_\varepsilon}{\varepsilon} &= 0, \\ \lim_{\varepsilon \to 0} \frac{x_\varepsilon - y_\varepsilon}{\varepsilon} &= 0, \\ \lim_{\varepsilon \to 0} \frac{\varepsilon^3}{(x_\varepsilon - y_\varepsilon)^2 + \varepsilon^2} &= 0.\end{aligned} \tag{2.9}$$

To show this we utilize the inequality

$$\begin{aligned}\Phi_\varepsilon(x_\varepsilon, t_\varepsilon, y_\varepsilon, s_\varepsilon) &\geq \Phi_\varepsilon(x_\varepsilon, t_\varepsilon, x_\varepsilon, t_\varepsilon) \\ &= u(x_\varepsilon, t_\varepsilon) - v(x_\varepsilon, t_\varepsilon) - 2\lambda t_\varepsilon - 2\varepsilon x_\varepsilon^2 - \frac{2}{\varepsilon}\frac{1}{x_\varepsilon^2} - \varepsilon,\end{aligned}$$

which yields

$$\frac{(t_\varepsilon - s_\varepsilon)^2 + (x_\varepsilon - y_\varepsilon)^2}{\varepsilon^2} + \frac{\varepsilon^3}{(x_\varepsilon - y_\varepsilon)^2 + \varepsilon^2} \leq v(x_\varepsilon, t_\varepsilon) - v(y_\varepsilon, s_\varepsilon) + \lambda(t_\varepsilon - s_\varepsilon) + \varepsilon\left(x_\varepsilon^2 - y_\varepsilon^2\right) + \frac{1}{\varepsilon}\left(\frac{1}{x_\varepsilon^2} - \frac{1}{y_\varepsilon^2}\right) + \varepsilon.$$

Since the viscosity solution $v$ is uniformly continuous, this is easily seen to imply (2.9). Now we define the sets

$$U_\varepsilon = \left\langle \frac{x_\varepsilon}{2}, \frac{3x_\varepsilon}{2} \right\rangle \times \left\langle \frac{t_\varepsilon}{2}, \frac{3t_\varepsilon}{2} \right\rangle, \qquad V_\varepsilon = \left\langle \frac{y_\varepsilon}{2}, \frac{3y_\varepsilon}{2} \right\rangle \times \left\langle \frac{s_\varepsilon}{2}, \frac{3s_\varepsilon}{2} \right\rangle,$$

$$\tilde{U}_\varepsilon = \left[ \frac{x_\varepsilon}{4}, \frac{7x_\varepsilon}{4} \right] \times \left[ \frac{t_\varepsilon}{4}, \frac{7t_\varepsilon}{4} \right], \qquad \tilde{V}_\varepsilon = \left[ \frac{y_\varepsilon}{4}, \frac{7y_\varepsilon}{4} \right] \times \left[ \frac{s_\varepsilon}{4}, \frac{7s_\varepsilon}{4} \right].$$

Let $\chi_{1,\varepsilon}$ and $\chi_{2,\varepsilon}$ be two smooth functions such that

$$\chi_{1,\varepsilon}(x,t) = \begin{cases} 1 & \text{if } (x,t) \in U_\varepsilon, \\ 0 & \text{if } (x,t) \notin \tilde{U}_\varepsilon, \end{cases} \qquad \chi_{2,\varepsilon}(y,s) = \begin{cases} 1 & \text{if } (y,s) \in V_\varepsilon, \\ 0 & \text{if } (y,s) \notin \tilde{V}_\varepsilon, \end{cases}$$

and $0 \leq \chi_{(1,2),\varepsilon} \leq 1$. Using $\chi_{1,\varepsilon}$ and $\chi_{2,\varepsilon}$ we define the test functions

$$\varphi_\varepsilon(x,t) = \chi_{1,\varepsilon}(x,t) \Bigg( v(y_\varepsilon, s_\varepsilon) + \lambda(t + s_\varepsilon) + \varepsilon\left(x^2 + y_\varepsilon^2\right) + \frac{(t - s_\varepsilon)^2 + (x - y_\varepsilon)^2}{\varepsilon^2} + \frac{1}{\varepsilon}\left(\frac{1}{x^2} + \frac{1}{y_\varepsilon^2}\right) + \frac{\varepsilon^3}{(x - y_\varepsilon)^2 + \varepsilon^2} \Bigg), \tag{2.10}$$

$$\psi_\varepsilon(y,s) = \chi_{2,\varepsilon}(y,s) \Bigg( u(x_\varepsilon, t_\varepsilon) - \lambda(t_\varepsilon + s) - \varepsilon\left(x_\varepsilon^2 + y^2\right) - \frac{(t_\varepsilon - s)^2 - (x - y_\varepsilon)^2}{\varepsilon^2} - \frac{1}{\varepsilon}\left(\frac{1}{x_\varepsilon^2} + \frac{1}{y^2}\right) - \frac{\varepsilon^3}{(x_\varepsilon - y)^2 + \varepsilon^2} \Bigg). \tag{2.11}$$

Now, since $u - \varphi$ has a local maximum in $(x_\varepsilon, t_\varepsilon)$ and $v - \psi$ has a local minimum in $(y_\varepsilon, s_\varepsilon)$, and $a$ is continuous in $x_\varepsilon$ and $y_\varepsilon$ we have that

$$\varphi_{\varepsilon,t}(x_\varepsilon, t_\varepsilon) + H(\varphi_{\varepsilon,x}(x_\varepsilon, t_\varepsilon), a(x_\varepsilon), g(t_\varepsilon)) \leq 0,$$

and

$$\psi_{\varepsilon,s}(y_\varepsilon, s_\varepsilon) + H(\psi_{\varepsilon,y}(y_\varepsilon, s_\varepsilon), a(y_\varepsilon), g(s_\varepsilon)) \geq 0.$$

Subtracting these, and using first the Lipschitz continuity of $H$ (c.f. (**A**.**1**)), and then the Lipschitz continuity of $a$ and $g$ (c.f. (**A**.**6**), (**A**.**7**)) and the fact that $t_\varepsilon$ and $s_\varepsilon$ both are less than $S < t_1$, we find

$$\begin{aligned} \varphi_{\varepsilon,t}(x_\varepsilon, t_\varepsilon) - \psi_{\varepsilon,s}(y_\varepsilon, s_\varepsilon) &\leq H(\varphi_{\varepsilon,x}(x_\varepsilon, t_\varepsilon), a(x_\varepsilon, g(t_\varepsilon))) - H(\psi_{\varepsilon,y}(y_\varepsilon, s_\varepsilon), a(y_\varepsilon), g(s_\varepsilon)) \\ &\leq C\Big( |\varphi_\varepsilon(x_\varepsilon, t_\varepsilon) - \psi_{\varepsilon,y}(y_\varepsilon, s_\varepsilon)| \\ &\qquad + (|g(t_\varepsilon) - g(s_\varepsilon)| + |a(x_\varepsilon) - a(y_\varepsilon)|)(|\varphi_{\varepsilon,x}(x_\varepsilon, t_\varepsilon)| + |\psi_{\varepsilon,y}(y_\varepsilon, s_\varepsilon)|) \Big) \\ &\leq C\Big( |\varphi_\varepsilon(x_\varepsilon, t_\varepsilon) - \psi_{\varepsilon,y}(y_\varepsilon, s_\varepsilon)| \\ &\qquad + (|t_\varepsilon - s_\varepsilon| + |x_\varepsilon - y_\varepsilon|)(|\varphi_{\varepsilon,x}(x_\varepsilon, t_\varepsilon)| + |\psi_{\varepsilon,y}(y_\varepsilon, s_\varepsilon)|) \Big). \end{aligned}$$

We also have that

$$\begin{aligned}
\varphi_{\varepsilon,t}\left(x_\varepsilon,t_\varepsilon\right) &= \lambda + 2\frac{t_\varepsilon - s_\varepsilon}{\varepsilon^2},\\
\psi_{\varepsilon,s}\left(y_\varepsilon,s_\varepsilon\right) &= -\lambda + 2\frac{t_\varepsilon - s_\varepsilon}{\varepsilon^2},\\
\varphi_{\varepsilon,x}\left(x_\varepsilon,t_\varepsilon\right) &= 2\varepsilon x_\varepsilon + 2\frac{x_\varepsilon - y_\varepsilon}{\varepsilon^2} - \frac{2}{\varepsilon x_\varepsilon^3} - 2\frac{\varepsilon^3\left(x_\varepsilon - y_\varepsilon\right)}{\left(\left(x_\varepsilon - y_\varepsilon\right)^2 + \varepsilon^2\right)^2},\\
\psi_{\varepsilon,y}\left(y_\varepsilon,s_\varepsilon\right) &= -2\varepsilon y_\varepsilon + 2\frac{x_\varepsilon - y_\varepsilon}{\varepsilon^2} + \frac{2}{\varepsilon x_\varepsilon^3} - 2\frac{\varepsilon^3\left(x_\varepsilon - y_\varepsilon\right)}{\left(\left(x_\varepsilon - y_\varepsilon\right)^2 + \varepsilon^2\right)^2},
\end{aligned}$$

and so

$$\varphi_{\varepsilon,t}\left(x_\varepsilon,t_\varepsilon\right) - \psi_{\varepsilon,s}\left(y_\varepsilon,s_\varepsilon\right) = 2\lambda, \tag{2.12}$$

$$\left|\varphi_{\varepsilon}\left(x_\varepsilon,t_\varepsilon\right) - \psi_{\varepsilon,y}\left(y_\varepsilon,s_\varepsilon\right)\right| = 2\varepsilon\left(x_\varepsilon + y_\varepsilon\right) + \frac{2}{\varepsilon}\left(\frac{1}{x_\varepsilon^3} + \frac{1}{y_\varepsilon^3}\right), \tag{2.13}$$

and

$$\begin{aligned}
\left|\varphi_{\varepsilon,x}\left(x_\varepsilon,t_\varepsilon\right)\right| + \left|\psi_{\varepsilon,y}\left(y_\varepsilon,s_\varepsilon\right)\right| &\le 2\varepsilon\left(x_\varepsilon + y_\varepsilon\right) + \frac{4\left|x_\varepsilon - y_\varepsilon\right|}{\varepsilon^2}\\
&\qquad + \frac{2}{\varepsilon}\left(\frac{1}{x_\varepsilon^3} + \frac{1}{y_\varepsilon^3}\right) + 4\frac{\varepsilon^3\left(x_\varepsilon - y_\varepsilon\right)}{\left(\left(x_\varepsilon - y_\varepsilon\right)^2 + \varepsilon^2\right)^2}.
\end{aligned} \tag{2.14}$$

Now we have that

$$\frac{1}{x_\varepsilon^3} \le \frac{\varepsilon^{3/2}}{C_1^3}, \quad \text{and} \quad \frac{1}{y_\varepsilon^3} \le \frac{\varepsilon^{3/2}}{C_1^3},$$

so the second term in (2.13) tends to zero as $\varepsilon \to 0$, and since $x_\varepsilon\ y_\varepsilon \le C/\sqrt{\varepsilon}$, the first term in (2.13) also tends to zero. Thus we have

$$\lim_{\varepsilon\to 0}\left|\varphi_{\varepsilon}\left(x_\varepsilon,t_\varepsilon\right) - \psi_{\varepsilon,y}\left(y_\varepsilon,s_\varepsilon\right)\right| = 0. \tag{2.15}$$

Furthermore, since

$$\frac{\varepsilon^3\left(x_\varepsilon - y_\varepsilon\right)}{\left(\left(x_\varepsilon - y_\varepsilon\right)^2 + \varepsilon^2\right)^2} \le \frac{\left|x_\varepsilon - y_\varepsilon\right|}{\varepsilon}\,\frac{\varepsilon^4}{\varepsilon^4} \to 0$$

as $\varepsilon \to 0$, also the last term (2.14) will tend to zero as $\varepsilon \to 0$. Recalling that $\left|t_\varepsilon - s_\varepsilon\right|/\varepsilon$ and $\left|x_\varepsilon - y_\varepsilon\right|/\varepsilon$ both tend to zero as $\varepsilon \to 0$, it follows that

$$\lim_{\varepsilon\to 0}\left(\left|t_\varepsilon - s_\varepsilon\right| + \left|x_\varepsilon - y_\varepsilon\right|\right)\left(\left|\varphi_{\varepsilon,x}\left(x_\varepsilon,t_\varepsilon\right)\right| + \left|\psi_{\varepsilon,y}\left(y_\varepsilon,s_\varepsilon\right)\right|\right) = 0, \tag{2.16}$$

Now we have established

$$\begin{aligned}
2\lambda \le \lim_{\varepsilon\to 0} C\Big(&\left|\varphi_{\varepsilon}\left(x_\varepsilon,t_\varepsilon\right) - \psi_{\varepsilon,y}\left(y_\varepsilon,s_\varepsilon\right)\right|\\
&+ \left(\left|t_\varepsilon - s_\varepsilon\right| + \left|x_\varepsilon - y_\varepsilon\right|\right)\left(\left|\varphi_{\varepsilon,x}\left(x_\varepsilon,t_\varepsilon\right)\right| + \left|\psi_{\varepsilon,y}\left(y_\varepsilon,s_\varepsilon\right)\right|\right)\Big)\\
= 0,&
\end{aligned}$$

which is a contradiction, since $\lambda > 0$. This concludes the proof of Lemma 2.3. □

*Proof (of Theorem 2.1).* If $T \le t_1$, the statement of the theorem is identical to that of Lemma 2.3. If $T > t_1$, then by Lemma 2.3 $u(\cdot,t) \le v(\cdot,t)$ for $t \in [0,t_1]$. In the strip $[t_1,t_2]$ we can use $u(\cdot,t_1)$ and $v(\cdot,t_1)$ as initial data, and conclude that $u(\cdot,t) \le v(\cdot,t)$ also in the strip $[t_1,t_2]$. We can repeat this argument in each strip $[t_i,t_{i+1}]$ to prove the theorem for any finite $T$. □

By modifying the proof of Lemma 2.3, it is now straightforward to show stability, namely

**Theorem 2.4 (Stability).** *Assume that (**A.1**), (**A.6**), (**A.7**) and (**A.8**) hold. If $u$ and $v$ are two viscosity solutions with initial data $u_0$ and $v_0$ (in $BUC(\mathbb{R})$) respectively, then*

$$\|u(\cdot,t)-v(\cdot,t)\|_{L^\infty(\mathbb{R})} \le \|u_0-v_0\|_{L^\infty(\mathbb{R})}\,. \tag{2.17}$$

*Proof.* The proof of this theorem consists in repeating the arguments in the proof of Lemma 2.3, but starting by using the penalization

$$\begin{aligned}\Phi_\varepsilon(x,t,y,s) &= u(x,t)-v(y,s)-\lambda(t+s)-\varepsilon(x^2+y^2)-\frac{(t-s)^2+(|x|-|y|)^2}{\varepsilon^2}\\ &\quad -\frac{1}{\varepsilon}\left(\frac{1}{x^2}+\frac{1}{y^2}\right)-\frac{\varepsilon^3}{(|x|-|y|)^2+\varepsilon^2}-K\chi(x,y)-\|u_0-v_0\|_{L^\infty(\mathbb{R})}\,,\end{aligned}$$

instead of (2.6). □

## 3. The front tracking scheme

In this section we show existence of a viscosity solution of (1.1). This is done by first considering a front tracking scheme, which yields approximate solutions to both (1.1) and (1.14). As all front tracking schemes, this one is based on the solution of Riemann problems, therefore we start by detailing this.

The Riemann problem for the conservation law is the initial value problem where $g$ is constant and $a$ and $p_0$ take two values, i.e.,

$$\begin{cases} p_t + H(p,a_l)_x = 0, & p(x,0)=p_l \quad \text{if } x<0,\\ p_t + H(p,a_r)_x = 0, & p(x,0)=p_r \quad \text{if } x\ge 0,\end{cases} \tag{3.1}$$

where $p_{l,r}$ and $a_{l,r}$ are constants. Since $g$ is constant, we have omitted the $g$ dependence of $H$ in our notation. The entropy solution to this problem is found by finding two $p$-values $p'_{l,r}$ such that the (scalar) Riemann problem with a flux function $H(p,a_l)$, a left state $p_l$ and right state $p'_l$ is solved by using waves of non-positive speed, and the Riemann problem with flux $H(p,a_r)$ and left state $p'_r$ and right state $p_r$ are solved using waves of non-negative speed only. Since the mapping $p \mapsto H(p,a)$ has a global maximum at $p=0$ for all $a$ and is even, these states are found as follows:

- **If** $p_l \le 0$**:** then $p'_l$ is in the set $[-p_l,\infty\rangle$.
- **If** $p_l > 0$**:** then $p'_l$ is in the set $[0,\infty\rangle$.
- **If** $p_r < 0$**:** then $\tilde{p_r}$ is in the set $\langle-\infty,0]$.
- **If** $p_r \ge 0$**:** then $p'_r$ is in the set $\langle-\infty,-p_r]$.

Furthermore, the Rankine–Hugoniot condition implies that

$$H\left(p'_l,a_l\right) = H\left(p'_r,a_r\right). \tag{3.2}$$

This is still not enough to give a unique solution, and in [7] the unique entropy solution is determined by the (unique) pair $(p'_l,p'_r)$ such that

$$|p'_l-p'_r|$$

is minimal. For the flux functions considered in this paper, we can always find a unique solution to the Riemann problem in this way. This solution will consist of $p$-waves, over which $a$ is constant, and the discontinuity in $a$, which we call an $a$-wave.

Although the solution of the Riemann problem will in general not be a monotone function of $x/t$, the flux, $H(p(x,t),a(x))$ will be monotone between the two values $H(p_l,a_l)$ and $H(p_r,a_r)$. This observation can be used to bound the solution of the Riemann problem. Let $G^\pm(h,a)$ be the two local inverses of $H$, i.e.,

$$G^+\left(H(p,a),a\right) = |p|\,, \quad \text{and} \quad G^-\left(H(p,a),a\right) = -|p|\,.$$

Since

$$\min\left\{H\left(p_l,a_l\right),H\left(p_r,a_r\right)\right\} \le H(p(x,t),a(x)) \le \max\left\{H\left(p_l,a_l\right),H\left(p_r,a_r\right)\right\},$$

we have that

$$\begin{aligned}\min\left\{G^-\left(H\left(p_l,a_l\right),a_r\right),G^-\left(H\left(p_r,a_r\right),a_l\right),-\left|p_l\right|,-\left|p_r\right|\right\}\\ \leq p(x,t)\leq \max\left\{G^+\left(H\left(p_l,a_l\right),a_r\right),G^+\left(H\left(p_r,a_r\right),a_l\right),\left|p_l\right|,\left|p_r\right|\right\}.\end{aligned} \tag{3.3}$$

By the assumption (**A**.**5**) we know that the values

$$\min\left\{G^-\left(H\left(p_l,a_l\right),a_r\right),G^-\left(H\left(p_r,a_r\right),a_l\right)\right\} \quad\text{and}\quad \max\left\{G^+\left(H\left(p_l,a_l\right),a_r\right),G^+\left(H\left(p_r,a_r\right),a_l\right)\right\}$$

are finite and bounded. By the special form of $H$, in particular (**A**.**3**), (**A**.**2**) and (**A**.**4**), we have the coarser (but simpler) bound

$$|p(x,t)|\leq G^+\left(\min\left\{H\left(p_l,a_l\right),H\left(p_r,a_r\right)\right\},\max\left\{a_l,a_r\right\}\right). \tag{3.4}$$

3.1. **Front tracking with constant $g$.** We start by defining the front tracking scheme for the case where $g$ is constant, this is a variation of the front tracking schemes defined in [13, 3]. Therefore consider the initial value problem

$$\begin{cases}p_t+H(p,a)_x=0 & \text{for } x\in\mathbb{R},\ t>0,\\ p(x,0)=p_0(x) & \text{for } x\in\mathbb{R}.\end{cases} \tag{3.5}$$

Let

$$z(p,a)=-\operatorname{sign}(p)\left(H(p,a)-H(0,a)\right) \quad\text{and}\quad \alpha(a)=H(0,a).$$

Since $a\mapsto H(0,a)$ is non-decreasing, $a\mapsto\alpha(a)$ is invertible. In the $(z,\alpha)$ plane, $a$ waves are straight lines of slope $\pm1$. An $a$-wave connecting two points $(z_1,\alpha_1)$ and $(z_2,\alpha_2)$ have slope 1 if $z_1$ and $z_2$ are non-positive, and slope $-1$ if these values are non-negative. If $z_1$ and $z_2$ have different sign, there is no $a$-wave connecting these points. Since $p$-waves connect points with the same $a$ values, these are horizontal lines in the $(z,\alpha)$ plane. Now fix a (small) number $\delta>0$, and set $\alpha_i=i\delta$, and $z_j=j\delta$, for integers $i$ and $j$. We define $p_0^\delta$ and $a^\delta$ as piecewise constant functions, with a finite number of jump discontinininuities, such that

$$\left.\begin{aligned}\left\|a-a^\delta\right\|_{L^1(\mathbb{R})}\to 0,\\ \left\|p_0-p_0^\delta\right\|_{L^1(\mathbb{R})}\to 0\end{aligned}\right\}\quad \text{as } \delta\to 0. \tag{3.6}$$

Label the (finite number of) values of $p^\delta$ and $a^\delta$ $p_1,\dots,p_M$, and $a_1,\dots,a_N$ respectively. Let $\alpha_j$ be the $j$th member of the ordered set

$$\left\{\alpha_k\right\}_{k=m'}^{M'}\cup\left\{\alpha(a_k)\right\}_{k=1}^{M},$$

where $m'$ and $M'$ are chosen such that

$$m'\leq\min_x\alpha(a^\delta(x))<\max_x\alpha(a^\delta(x))\leq M'.$$

For ease of notation, set

$$a_j=\alpha^{-1}\left(\alpha_j\right).$$

Next for each $\alpha_j$, we define $z_{j,k}$ to be the $k$th member of the ordered set

$$\left\{z_i\right\}_{i=-N'(j)}^{N'(j)}\cup\left\{z\left(p_i,a_j\right)\right\}_{i=1}^{M},$$

where $N'(j)$ is such that

$$z^{-1}\left(z_{-N'(j)},a_j\right)=-P, \quad\text{and}\quad z^{-1}\left(z_{N'(j)},a_j\right)=P,$$

where the value $P$ will be determined below. We also set

$$p_{j,k}=z^{-1}\left(z_{j,k},a_j\right), \quad\text{and}\quad H_{j,k}=H\left(p_{j,k},a_j\right).$$

Then, for each $j$, let the approximate flux function $H^\delta(p,a)$ be the piecewise linear interpolant,

$$H^\delta\left(p,a_j\right)=H_{j,k}+\left(p-p_{j,k}\right)\frac{H_{j,k+1}-H_{j,k}}{p_{j,k+1}-p_{j,k}}, \quad \text{for } p\in[p_{j,k},p_{j,k+1}]. \tag{3.7}$$

Now the front tracking solution, which we label $p^\delta=p^\delta(x,t)$ is constructed as follows. At $t=0$ we solve the Riemann problems defined by the discontinuities in $a^\delta$ and $p_0^\delta$, using the flux function $H^\delta$. The flux function $H^\delta$ is constructed so that $p$-rarefaction waves are now a series of contact discontinuities. Thus the solution of the Riemann problem is a piecewise constant function of $x/t$. Furthermore, all the intermediate states will be on the grid $(p_{j,k}, a_j)$. At some time $t > 0$, two discontinuities, hereafter called *fronts*, will collide, giving a new Riemann problem centered at the interaction points. We solve this Riemann problem, this gives new fronts that can be tracked until the next interaction point and so on. Based on the estimate (3.4), we also have that

$$|p^\delta(x,t)| \le G^+ \left(\inf_x H\left(p_0^\delta(x), a^\delta(x)\right), \sup a^\delta(x)\right) =: P. \tag{3.8}$$

In order to show that $p^\delta$ is well-defined (for instance we must show that there is no accumulation of collision times) we define the Temple functional of a front by

$$T(w) = \begin{cases} |\Delta\Psi| & \text{if } w \text{ is a } u\text{-front}, \\ 2\,|\Delta a|\, g & \text{if } w \text{ is an } a\text{-front, and } \Psi_r < \Psi_l, \\ 4\,|\Delta a|\, g & \text{if } w \text{ is an } a\text{-front, and } \Psi_r > \Psi_l, \end{cases} \tag{3.9}$$

where $\Psi$ is defined in (1.15), and we have included the $g$ value in our notation since we shall need it later. For sequence of fronts, define $T$ additively. Next, for the front tracking approximation $p^\delta$, $T$ is defined as the sum over all the fronts in $p^\delta$, and with a slight abuse of notation we write $T(p^\delta)$. From estimates found in [20, 13] it follows that $t \mapsto T(p^\delta(\cdot,t))$ is non-increasing, which again implies that $\left|\Psi(p^\delta, a^\delta)\right|_{BV}$ is bounded. This again implies that $p^\delta$ is well defined, and furthermore that there is only a finite number of interactions of fronts for all $t > 0$.

Summing up, we have chosen the grid so that the entropy solution to the initial value problem

$$\begin{aligned} p_t + H^\delta\left(p, a^\delta\right)_x = 0, \quad t > 0, \quad x \in \mathbb{R} \\ p(x,0) = p_0^\delta(x), \quad x \in \mathbb{R}, \end{aligned} \tag{3.10}$$

can be constructed by front tracking for any time $t$. Furthermore $p^\delta$ will take values that are grid points, i.e., for any point $(x,t)$ such that $p^\delta$ and $a^\delta$ are constant at $(x,t)$,

$$z\left(p^\delta\left(x,t\right), a^\delta(x)\right) = z_{j,k}, \qquad \text{for some } j \text{ and } k.$$

In particular, this means that

$$H^\delta\left(p^\delta, a^\delta\right) = H\left(p^\delta, a^\delta\right), \quad \text{almost everywhere}.$$

For an elaboration and proof of these statements, see [13]. The construction used here differs from the construction in [13] in that we have added grid points corresponding to the discretization of the initial function $p_0$ and the coefficient $a$, instead of choosing discretization that take values on the fixed grid in the $(z, \alpha)$ plane.

3.2. **Front tracking in general.** Now we can define the front tracking approximation in the case where $g$ is not constant, c.f. (1.14). Let $g^\delta$ be a piecewise constant approximation to $g$, such that

$$\begin{aligned} \left\|g^\delta - g\right\|_{L^1(\mathbb{R}^+)} \to 0, \quad \text{as } \delta \to 0, \\ \left|g^\delta\right|_{BV(\langle 0,T])} \le |g|_{BV(\langle 0,T])}. \end{aligned} \tag{3.11}$$

Define $t^n$ such that $g^\delta$ is constant on each interval $I^n = \langle t^n, t^{n+1}]$. Assuming that we can define front tracking for $t < t^n$, we can then use $p^\delta(\cdot, t^n)$ as initial values for a front tracking approximation defined in $[t^n, t^{n+1}\rangle$. In order to do this we must use a "new" mapping $z$, since $z = z(p, a, g)$, and redefine the grid on which we operate. However, we keep the grid points corresponding to $p^\delta(\cdot, t^n)$. In this way, the grid used in the interval $I^{n+1}$ will contain more points than the one used in $I^n$, but since there are only a finite number of intervals $I^n$ such that $t^n \le T$, for a fixed $\delta$, we use a finite number of grid points for $t \le T$. If, for $t \in I^n$, $H^\delta(\cdot, \cdot, g^\delta(t))$ denotes the approximate flux function constructed above using $H(\cdot, \cdot, g^\delta\,|_{I^n})$ and $p^\delta(\cdot, t^n)$, then we have that the front tracking construction $p^\delta$ will be an entropy solution of

$$\begin{aligned} p_t^\delta + H^\delta\left(p^\delta, a^\delta(x), g^\delta(t)\right)_x = 0, \quad t > 0, \quad x \in \mathbb{R} \\ p^\delta(x,0) = p_0^\delta(x), \quad x \in \mathbb{R}. \end{aligned} \tag{3.12}$$

We call the discontinuities in $u^\delta$ *fronts*, and we have three types, $u$-fronts, $a$-fronts and $g$-fronts (that have infinite speed!).

3.3. **Compactness.** We aim to show that the sequence $\{p^\delta\}_{\delta>0}$ is compact in $L^1$, by estimating the variation of $\Psi(p^\delta, a^\delta, g^\delta)$. For each time $t$, such that $g^\delta$ is constant at $t$, we can view $p^\delta$ as consisting of a sequence of fronts, $u$-fronts and $a$-fronts, and we know that $T(p^\delta)$ is non-increasing in $I^n$.

For $t \in I^n$ we also know that

$$H\left(p^\delta(x,t), a^\delta(x), g^n\right) \geq \inf_x H\left(p^\delta\left(x,t_n\right), a^\delta(x), g^n\right).$$

Now we define the sequence $\{\bar{p}^n\}_{n\geq 1}$ by

$$\begin{cases} \bar{p}^1 = G^+\left(\inf_x H\left(p_0^\delta(x), a^\delta(x), g^1\right)\right), \bar{a}, g^1\right) \\ \bar{p}^k = G^+\left(H\left(\bar{p}^{k-1}, \underline{a}, g^k\right)\right), \bar{a}, g^k\right) & \text{for } k > 1, \end{cases}$$

where $\underline{a} = \inf_x a(x)$, and $\bar{a} = \sup_x a(x)$, and $g^\delta(t) = g^n$ for $t \in I^n$. By the continuity of $H$,

$$H\left(\bar{p}^n, \underline{a}, g^{n+1}\right) \geq H\left(\bar{p}^n, \underline{a}, g^n\right) - C\left|\Delta g\right|,$$

where $\Delta g = g^{n+1} - g^n$. Also

$$G^+(h, a, g^{n+1}) \leq G^+(h, a, g^n) + C\left|\Delta g\right|,$$

for some constant $C$. This means that

$$\bar{p}^{n+1} \leq \bar{p}^n + C\left|\Delta g\right|, \quad \text{and thus} \quad \bar{p}^n \leq \bar{p}^1 + C\left|g\right|_{BV}.$$

Now, for $t \in I^n$,

$$\begin{aligned} \left|p^\delta(x,t)\right| &\leq G^+\left(\inf_x H\left(p^\delta(x,t_n), a^\delta(x), g^n\right), \bar{a}, g^n\right) \\ &\leq G^+\left(H\left(\bar{p}^{n-1}, \underline{a}, g^n\right), \bar{a}, g^n\right) \\ &= \bar{p}^n. \end{aligned}$$

Hence the sequence of front tracking approximations $\{p^\delta\}$ is uniformly bounded independently of $\delta$, i.e.,

$$\left|p^\delta(x,t)\right| \leq \bar{C}\ (= \bar{C}(p_0, H, a, g)). \tag{3.13}$$

where $\bar{C}(\cdots)$ is a constant (depending on its arguments only).

Now we define a "Glimm type" functional, which we shall show that is nonincreasing in time, and then this will imply that the total variation of $\Psi(p^\delta, a, g)$ is bounded. Set

$$Q(t) = T(t)\left|g^\delta(\cdot)\right|_{BV([t,T])}, \tag{3.14}$$

where with a slight abuse of notation we write $T(t) = T(p^\delta(\cdot, t))$. With these definitions, we can state the following lemma.

**Lemma 3.1.** *There exists a positive constant $C$, depending only on $H$, $a$ and $g$, such that for all $t > 0$, we have that the "Glimm functional"*

$$G(t) = T(t) + CQ(t), \tag{3.15}$$

*is nonincreasing in time.*

*Proof.* The proof of this lemma is very similar to the proof of the corresponding lemma in [3], and we detail only the differences.

In each interval $I^n$, we know from [13] that $T$ is non-increasing, and the lemma holds. To prove the lemma we must study interactions between $p$-fronts and $g$-fronts, and between $p$-fronts and $g$-fronts.

Now the proof of the lemma for the interaction of a $p$-front and an $a$-front is identical to the proof of the corresponding case in [3], which means that we only must study the interaction of a $p$-front with a $g$-front.

First note that by (**A.3**) and (**A.1**) and (**A.4**) there exists a value $P$ and positive constants $c_0$ and $C_0$ such that

$$\begin{array}{rl} |p| \leq P & \Rightarrow 0 < c_0 \leq H_{pp}(p,a,g) \leq C_0 \\ P \leq |p| \leq \bar{C} & \Rightarrow 0 < c_0 \leq |H_p(p,a,g)| \leq C_0 \end{array} \quad \text{if } \underline{a} \leq a \leq \bar{a} \text{ and } \underline{g} \leq g \leq \bar{g}, \tag{3.16}$$

where $\underline{g} = \inf g$ and $\bar{g} = \sup g$. Now we consider the interaction of a single $u$-wave and a single $g$-wave. The situation is depicted in Figure 1. For this interaction we claim that

$$\begin{aligned} \left|\Psi\left(p_r,a,g^+\right)-\Psi\left(p_l,a,g^+\right)\right| &- \left|\Psi\left(p_r,a,g^-\right)-\Psi\left(p_l,a,g^-\right)\right| \\ &\leq C\left|g^+-g^-\right|\left|\Psi\left(p_r,a,g^-\right)-\Psi\left(p_l,a,g^-\right)\right|. \end{aligned} \tag{3.17}$$

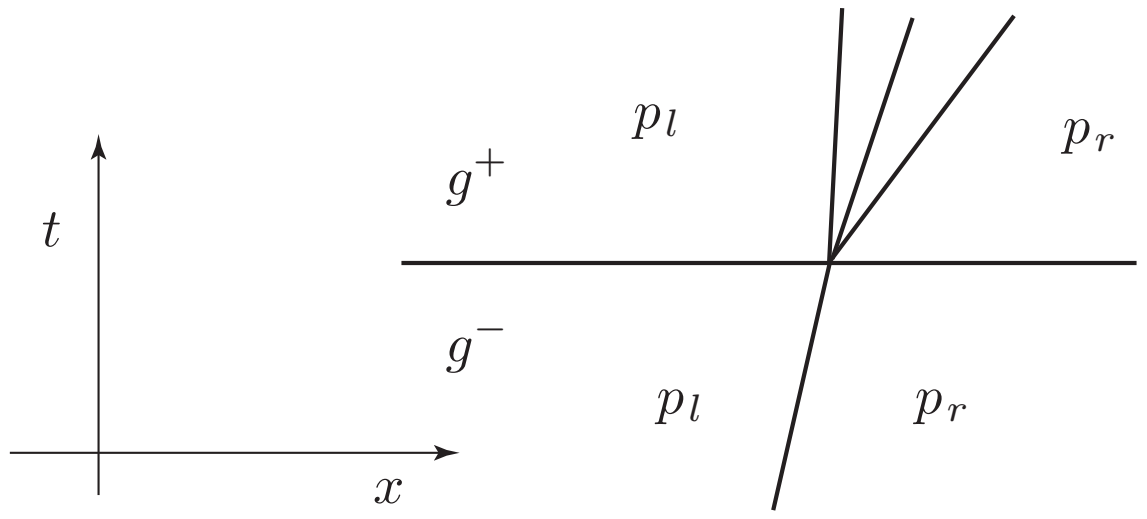

FIGURE 1. The states used in an interaction between a $p$-wave and a $g$-wave

We start by noting that $\Psi(0,\cdot,\cdot) = \Psi_p(0,\cdot,\cdot) = 0$, and

$$\begin{aligned} \left|\Psi\left(p_r,a,g^+\right)-\Psi\left(p_l,a,g^+\right)\right| &- \left|\Psi\left(p_r,a,g^-\right)-\Psi\left(p_l,a,g^-\right)\right| \\ &\leq \left|\Psi\left(p_r,a,g^+\right)-\Psi\left(p_l,a,g^+\right)-\Psi\left(p_r,a,g^-\right)-\Psi\left(p_l,a,g^-\right)\right|. \end{aligned}$$

To prove (3.17) we consider different cases.

**Case 1:** $|p_l| \leq P$ and $|p_r| \leq P$. Now

$$\begin{aligned} \Psi\left(p_r,a,g^+\right) &- \Psi\left(p_l,a,g^+\right)-\Psi\left(p_r,a,g^-\right)-\Psi\left(p_l,a,g^-\right) \\ &= \int_{p_l}^{p_r} \left(\Psi_p\left(\xi,a,g^+\right)-\Psi_p\left(\xi,a,g^-\right)\right)\,d\xi \\ &= \int_{p_l}^{p_r} \left(\Psi_p\left(\xi,a,g^+\right)-\Psi_p\left(0,a,g^+\right)-\Psi_p\left(\xi,a,g^-\right)+\Psi_p\left(0,a,g^-\right)\right)\,d\xi \\ &= \int_{p_l}^{p_r}\int_0^{\xi} \left(\Psi_{pp}\left(\eta,a,g^+\right)-\Psi_{pp}\left(\eta,a,g^-\right)\right)\,d\eta d\xi \\ &= \int_{p_l}^{p_r}\int_0^{\xi}\int_{g^-}^{g^+} \Psi_{ppg}(\eta,a,g)\,dgd\eta d\xi, \end{aligned}$$

and

$$\begin{aligned}\Psi\left(p_r,a,g^-\right)-\Psi\left(p_l,a,g^-\right)&=\int_{p_l}^{p_r}\Psi_p\left(\xi,a,g^-\right)\,d\xi\\&=\int_{p_l}^{p_r}\left(\Psi_p\left(\xi,a,g^-\right)-\Psi_p\left(0,a,g^-\right)\right)\,d\xi\\&=\int_{p_l}^{p_r}\int_0^{\xi}\Psi_{pp}\left(\eta,a,g^-\right)\,d\eta d\xi.\end{aligned}$$

Moreover observe that

$$\Psi_{pp}(p,a,g)=-\mathrm{sign}\,(p)\,\frac{H_{pp}(p,a,g)}{H(0,a,g)} \tag{3.18}$$

and

$$\Psi_{ppg}(p,a,g)=\mathrm{sign}\,(p)\,\frac{H_{pp}(p,a,g)H_g(0,a,g)-H_{ppg}(p,a,g)H(0,a,g)}{H(0,a,g)^2}. \tag{3.19}$$

To fix ideas, we assume that $p_l<p_r$, so that by (**A.2**), $\Psi(p_l,\cdot,\cdot)\le\Psi(p_r,\cdot,\cdot)$. If $0\le p_l<p_r$ then

$$\begin{aligned}&\left|\Psi\left(p_r,a,g^+\right)-\Psi\left(p_l,a,g^+\right)-\Psi\left(p_r,a,g^-\right)-\Psi\left(p_l,a,g^-\right)\right|\\&\qquad\le C_1\left|g^+-g^-\right|\int_{p_l}^{p_r}\int_0^{\xi}d\eta d\xi\\&\qquad=C_1\left|g^+-g^-\right|\frac{p_r^2-p_l^2}{2},\end{aligned}$$

and

$$\Psi\left(p_r,a,g^-\right)-\Psi\left(p_l,a,g^-\right)\ge c_1\int_{p_l}^{p_r}\int_0^{\xi}d\eta d\xi,$$

for some positive constants $c_1$ and $C_1$. Therefore

$$\begin{aligned}&\left|\Psi\left(p_r,a,g^+\right)-\Psi\left(p_l,a,g^+\right)-\Psi\left(p_r,a,g^-\right)-\Psi\left(p_l,a,g^-\right)\right|\\&\qquad\le\frac{C_1}{c_1}\left|g^+-g^-\right|\left|\Psi\left(p_r,a,g^-\right)-\Psi\left(p_l,a,g^-\right)\right|,\end{aligned} \tag{3.20}$$

and thus the claim holds. Next, if $p_l<p_r\le 0$, then

$$\begin{aligned}&\left|\Psi\left(p_r,a,g^+\right)-\Psi\left(p_l,a,g^+\right)-\Psi\left(p_r,a,g^-\right)-\Psi\left(p_l,a,g^-\right)\right|\\&\qquad\le C_2\left|g^+-g^-\right|\int_{p_l}^{p_r}\int_{\xi}^{0}d\eta d\xi,\end{aligned}$$

and

$$\Psi\left(p_r,a,g^-\right)-\Psi\left(p_l,a,g^-\right)\ge c_2\int_{p_l}^{p_r}\int_{\xi}^{0}d\eta d\xi,$$

for some positive constants $C_2$ and $c_2$. Hence (3.20) holds also in this case. If $p_l \leq 0 \leq p_r$ then we write

$$\begin{aligned}&\left|\Psi\left(p_r,a,g^+\right)-\Psi\left(p_l,a,g^+\right)-\Psi\left(p_r,a,g^-\right)-\Psi\left(p_l,a,g^-\right)\right|\\&\qquad\leq C_3\left|g^+-g^-\right|\left(\int_{p_l}^{0}\int_{\xi}^{0}d\eta d\xi+\int_{0}^{p_r}\int_{0}^{\xi}d\eta d\xi\right),\end{aligned}$$

and

$$\Psi\left(p_r,a,g^-\right)-\Psi\left(p_l,a,g^-\right)\geq c_3\left(\int_{p_l}^{0}\int_{\xi}^{0}d\eta d\xi+\int_{0}^{p_r}\int_{0}^{\xi}d\eta d\xi\right),$$

for some positive constants $C_3$ and $c_3$, and (3.20) follows. If $p_r < p_l$ we can use the same arguments.

**Case 2:** $|p_l| \leq P$ and $|p_r| \leq P$. We start by observing that

$$\Psi\left(p_r,a,g^+\right)-\Psi\left(p_l,a,g^+\right)-\Psi\left(p_r,a,g^-\right)-\Psi\left(p_l,a,g^-\right)=\int_{p_l}^{p_r}\int_{g^-}^{g^+}\Psi_{pg}(\xi,a,\theta)\,d\theta d\xi,$$

and

$$\Psi\left(p_r,a,g^-\right)-\Psi\left(p_l,a,g^-\right)=\int_{p_l}^{p_r}\Psi_p(\xi,a,g^-)\,d\xi.$$

Since $\Psi_{pg}$ is bounded, we have that

$$\left|\Psi\left(p_r,a,g^+\right)-\Psi\left(p_l,a,g^+\right)-\Psi\left(p_r,a,g^-\right)-\Psi\left(p_l,a,g^-\right)\right|\leq C_4\left|g^+-g^-\right|\left|p_r-p_l\right|,$$

and since $\Psi_p \geq c_4$ for $p \notin \langle -P, P\rangle$,

$$\left|\Psi\left(p_r,a,g^-\right)-\Psi\left(p_l,a,g^-\right)\right|\geq c_4\left|p_r-p_l\right|,$$

for some positive constants $c_4$ and $C_4$. Thus (3.20) follows.

**Case 3:** $|p_l| \leq P \leq p_r$. In this case we start by writing

$$\begin{aligned}&\Psi\left(p_r,a,g^+\right)-\Psi\left(p_l,a,g^+\right)-\Psi\left(p_r,a,g^-\right)+\Psi\left(p_l,a,g^-\right)\\&\qquad=\Psi\left(p_r,a,g^+\right)-\Psi\left(P,a,g^+\right)+\Psi\left(P,a,g^+\right)-\Psi\left(p_l,a,g^+\right)\\&\qquad\quad-\Psi\left(p_r,a,g^-\right)+\Psi\left(P,a,g^-\right)-\Psi\left(P,a,g^-\right)+\Psi\left(p_l,a,g^-\right)\\&\qquad=\int_{P}^{p_r}\int_{g^-}^{g^+}\Psi_{pg}(\xi,a,\theta)\,d\theta d\xi+\int_{p_l}^{P}\int_{0}^{\xi}\int_{g^-}^{g^+}\Psi_{ppg}(\eta,a,\theta)\,d\theta d\eta d\xi,\end{aligned}$$

and

$$\begin{aligned}\Psi\left(p_r,a,g^-\right)-\Psi\left(p_l,a,g^-\right)&=\Psi\left(p_r,a,g^-\right)-\Psi\left(P,a,g^-\right)+\Psi\left(P,a,g^-\right)-\Psi\left(p_l,a,g^-\right)\\&=\int_{P}^{p_r}\Psi_p\left(\xi,a,g^-\right)\,d\xi+\int_{p_l}^{P}\int_{0}^{\xi}\Psi_{pp}\left(\eta,a,g^-\right)\,d\eta d\xi.\end{aligned}$$

Since the derivatives of $\Psi$ are bounded,

$$\left|\Psi\left(p_r,a,g^+\right)-\Psi\left(p_l,a,g^+\right)-\Psi\left(p_r,a,g^-\right)+\Psi\left(p_l,a,g^-\right)\right|\leq C_5\left(p_r-P+\frac{P^2-p_l^2}{2}\right)\left|g^+-g^-\right|,$$

for some positive constant $C_5$. Since $\Psi_{pp}$ is strictly positive inside $[-P, P]$ and $\Psi_p$ is larger than some fixed constant outside this interval,

$$\left|\Psi\left(p_r,a,g^-\right)-\Psi\left(p_l,a,g^-\right)\right|\geq c_5\left(p_r-P+\frac{P^2-p_l^2}{2}\right),$$

for some positive constant $c_5$, and thus (3.20) holds. If $p_l < 0$, then

$$\left|\Psi\left(p_r,a,g^+\right)-\Psi\left(p_l,a,g^+\right)-\Psi\left(p_r,a,g^-\right)+\Psi\left(p_l,a,g^-\right)\right|\le C_6\left(p_r-P+\frac{P^2+p_l^2}{2}\right)\left|g^+-g^-\right|,$$

and

$$\left|\Psi\left(p_r,a,g^-\right)-\Psi\left(p_l,a,g^-\right)\right|\ge c_6\left(p_r-P+\frac{P^2+p_l^2}{2}\right),$$

so (3.20) holds again.

**Case 4:** $|p_l|\le P$ and $p_r\le -P$. This is analogous to case 3.

**Case 5:** $|p_r|\le P\le p_l$. This is analogous to case 3.

**Case 6:** $|p_r|\le P$ and $p_l\le -P$. This is analogous to case 3.

**Case 7:** $p_l\le -P$ and $p_r\ge P$. Now we write

$$\begin{aligned}
\Psi\left(p_r,a,g^+\right)&-\Psi\left(p_l,a,g^+\right)-\Psi\left(p_r,a,g^-\right)+\Psi\left(p_l,a,g^-\right)\\
&=\Psi\left(p_r,a,g^+\right)-\Psi\left(P,a,g^+\right)+\Psi\left(P,a,g^+\right)-\Psi\left(-P,a,g^+\right)\\
&\quad+\Psi\left(-P,a,g^+\right)-\Psi\left(p_l,a,g^+\right)-\Psi\left(p_r,a,g^-\right)+\Psi\left(P,a,g^-\right)\\
&\quad-\Psi\left(P,a,g^-\right)+\Psi\left(-P,a,g^-\right)-\Psi\left(-P,a,g^-\right)+\Psi\left(p_l,a,g^-\right)\\
&=\int_P^{p_r}\int_{g^-}^{g^+}\Psi_{pg}(\xi,a,\theta)\,d\theta d\xi+\int_{-P}^{P}\int_0^{\xi}\int_{g^-}^{g^+}\Psi_{ppg}(\eta,a,\theta)\,d\theta d\eta d\xi+\int_{p_l}^{-P}\int_{g^-}^{g^+}\Psi_{pg}(\xi,a,\theta)\,d\theta d\xi,
\end{aligned}$$

and

$$\begin{aligned}
0&\le\Psi\left(p_r,a,g^-\right)-\Psi\left(p_l,a,g^-\right)\\
&=\Psi\left(p_r,a,g^-\right)-\Psi\left(P,a,g^-\right)+\Psi\left(P,a,g^-\right)-\Psi\left(-P,a,g^-\right)\\
&\quad+\Psi\left(-P,a,g^-\right)-\Psi\left(p_l,a,g^-\right)\\
&=\int_P^{p_r}\Psi_p\left(\xi,a,g^-\right)\,d\xi+\int_{-P}^{P}\int_0^{\xi}\Psi_{pp}\left(\eta,a,g^-\right)\,d\eta d\xi+\int_{p_l}^{-P}\Psi_p\left(\xi,a,g^-\right)\,d\xi.
\end{aligned}$$

As in the earlier cases, (3.20) is straightforward to show from this.

**Case 8:** $p_r\le -P$ and $p_l\ge P$. This is analogous to case 7.

Now the proof of (3.17) and thereby of Lemma 3.1 is finished. □

Let $T^n=T\left|_{I^n}\right.$ and $g^n=g^\delta\left|_{I^n}\right.$. Since $T$ is non-increasing in each interval $I^n$, from Lemma 3.1, we have that

$$T^{n+1}\le T^n\left(1+C\left|g^{n+1}-g^n\right|\right).$$

By the Grönwall inequality it follows that

$$\begin{aligned}
T(t)&\le T^1(0+)\exp\left(\sum_n\left|g^n-g^{n-1}\right|\right)\\
&\le\lim_{s\downarrow 0}T(s)\exp\left(|g|_{BV}\right)\\
&\le\left(\left|\Psi\left(p_0,a,g(0)\right)\right|_{BV}+4\left|a\right|_{BV}\left|g(0)\right|\right)e^{|g|_{BV}}.
\end{aligned}\tag{3.21}$$

where the sum in the first line above is over those $n$ such that $t_n<t$.

This clearly implies that the total variation $\Psi(p^\delta,a^\delta,g^\delta(t))$ is bounded independently of $\delta$ and $t$. In particular, this means that the front tracking construction is well-defined, and we have a finite number of fronts and interaction of fronts, for $0\le t\le T$. For a proof of this, see e.g. [13]. Furthermore, since $\left|p^\delta\right|\le\bar{C}$, c.f. (3.13),

$$\Psi\left(\bar{C},\underline{a},\underline{g}\right)\le\Psi\left(u^\delta(x,t),a^\delta(x),g^\delta(t)\right)\le\Psi\left(\bar{C},\bar{a},\bar{g}\right).$$

By Helly's theorem, for each fixed $t\in[0,T]$,

$$\Psi\left(p^\delta(\cdot,t),a^\delta,g^\delta(t)\right)\to\psi,\quad\text{almost everywhere as }\delta\downarrow 0,$$

and by the Lebesgue's dominated convergence theorem also in $L^1(\mathbb{R})$. Furthermore, by a diagonal argument, we can achieve this convergence for a dense countable set $\{t^\gamma\} \subset [0,T]$. For $t^\gamma$ in this set, define

$$p(\cdot, t^\gamma) = \Psi^{-1}\left(\psi, a, g\left(t^\gamma\right)\right).$$

Hence also $p^\delta(\cdot, t^\gamma)$ converges to some $p(\cdot, t^\gamma)$. For any $t \in [0,T]$ we have that

$$\begin{aligned}\left\|p^{\delta_1}(\cdot,t) - p^{\delta_2}(\cdot,t)\right\|_{L^1(\mathbb{R})} &\leq \left\|p^{\delta_1}(\cdot,t^\gamma) - p^{\delta_1}(\cdot,t)\right\|_{L^1(\mathbb{R})} \\ &+ \left\|p^{\delta_1}(\cdot,t^\gamma) - p^{\delta_2}(\cdot,t^\gamma)\right\|_{L^1(\mathbb{R})} + \left\|p^{\delta_2}(\cdot,t^\gamma) - p^{\delta_2}(\cdot,t)\right\|_{L^1(\mathbb{R})},\end{aligned}$$

where $t^\gamma$ is such that $p^\delta(\cdot, t^\gamma) \to p(\cdot, t^\gamma)$. It is easy to show, as in [3], that the map $t \mapsto u^\delta(\cdot, t)$ is $L^1$ Lipschitz continuous, so the first and third terms above can be made arbitrarily small by choosing $\delta_1$ and $\delta_2$ small, and the middle term can be made small by choosing $t^\gamma$ close to $t$. Hence we have that $p^\delta$ converges to some function $p$ in $L^1(\mathbb{R} \times [0,T])$.

Now since $p \mapsto H(p,a,g)$ has a unique maximum for $p = 0$ for all $a$ and $g$, we can use the same arguments as in [3] to show that $p$ is an entropy solution to (1.14) in the sense of (1.17). Furthermore, using arguments from [11] and [3] the entropy solution is unique. Summing up, we have proved:

**Theorem 3.2.** *Assume that **(A.1)** – **(A.8)** all hold. Let $u_0 \in L^1(\mathbb{R})$ be such that the total variation of $\Psi(u_0, a, g)$ is bounded. Then there exists a unique entropy solution to* (1.14) *$p = p(x,t)$. This solution can be constructed as the limit of the front tracking scheme outlined above.*

3.4. **Front tracking for the Hamilton–Jacobi equation.** Now we show how the front tracking approximation to the entropy solution of the conservation law also yields a front tracking approximation to the viscosity solution of the Hamilton–Jacobi equation. We start by studying the Riemann problem.

**Lemma 3.3.** *Assume that **(A.5)**, **(A.4)**, **(A.3)** and **(A.1)** all hold, then the Riemann problem for the Hamilton–Jacobi equation*

$$\begin{cases} u_t + H\left(u_x, a_l\right) = 0 & \text{if } x \leq 0 \text{ and } t > 0, \\ u_t + H\left(u_x, a_r\right) = 0 & \text{if } x > 0 \text{ and } t > 0, \\ u(x,0) = u_0(0) + \begin{cases} p_l x & x \leq 0, \\ p_r x & x > 0, \end{cases} \end{cases} \tag{3.22}$$

*has a unique viscosity solution given by*

$$u(x,t) = u_0(0) + x p(x,t) - t \begin{cases} H(p(x,t), a_l) & x \leq 0 \\ H(p(x,t), a_r) & x > 0 \end{cases} \qquad \text{for } t > 0, \tag{3.23}$$

*where $p(x,t)$ is the unique entropy solution of the Riemann problem for the conservation law* (3.1)*.*

**Remark 3.4.** From (3.23) and the fact that $p$ solves a conservation law, it follows that

$$u(x,t) = u_0(0) - tH_0 + \int_0^x p(x,t)\,dx, \tag{3.24}$$

where $H_0 = H(p_l', a_l) = H(p_r', a_r)$, is an alternative formula for $u$. This can be shown by observing that $p = p(x/t)$, and differentiating (3.23) with respect to $x$, using that $p_t + H(p,a)_x = 0$.

*Proof.* Let $u$ be defined by (3.23). If $p$ is continuous at $(x,t)$ then so is $u$, and if $p$ has a discontinuity moving with speed $\sigma$, set $x = \sigma t$. By the Rankine–Hugoniot condition, if $\sigma \neq 0$,

$$\begin{aligned} u\left(x^-, t\right) &= u_0(0) + t\sigma p\left(x^-, t\right) - tH\left(p(x^-, t), a_{l,r}\right) \\ &= u_0(0) + t\left(\sigma p\left(x^-, t\right) - H\left(p(x^-, t), a_{l,r}\right)\right) \\ &= u_0(0) + t\left(\sigma p\left(x^+, t\right) - H\left(p(x^+, t), a_{l,r}\right)\right) = u\left(x^+, t\right), \end{aligned}$$

where we use $a_l$ if $\sigma < 0$ and $a_r$ if $\sigma > 0$. If $\sigma = 0$, then $H(p(0^-, t), a_l) = H(p(0^+, t), a_r)$ which gives

$$u\left(0^-, t\right) = u_0(0) + tH\left(p\left(0^-, t\right), a_l\right) = u_0(0) + tH\left(p\left(0^+, t\right), a_r\right) = u\left(0^+, r\right).$$

Thus, $u$ is uniformly continuous. Now let $\varphi$ be a test function, and assume that $u-\varphi$ has a maximum at $(x_0,t_0)$. We proceed by studying two cases.

**Case 1:** $x_0 \neq 0$. First assume that $x_0>0$, and let $q$ be a solution of the (scalar) Riemann problem,

$$q_t + H(q,a_r)_x = 0, \qquad q(x,0) = \begin{cases} p(t_0,0^+) & x \le 0, \\ p_r & x > 0. \end{cases}$$

Clearly $p(x,t)=q(x,t)$ for $x>0$ and $t>0$, So by [10, Proposition 2.3]

$$\varphi_t\left(x_0,t_0\right) + H\left(\varphi_x\left(x_0,t_0\right),a_r\right) \le 0.$$

If $x_0<0$, then we replace $q$ by the solution of the Riemann problem with initial data given by $p_l$ and $p(0^-,t_0)$ and flux $H(q,a_l)$, and reach the same conclusion.

**Case 2:** $x_0=0$. Now by (3.24) it follows that

$$p_l' = \lim_{x\uparrow 0} u_x(x,t_0), \quad \text{and} \quad p_r' = \lim_{x\downarrow 0} u_x(x,t_0),$$

where $p_{l,r}'$ are the states adjacent to $x=0$ of the Riemann solution $p$. Since $u-\varphi$ has a local maximum,

$$p_r' \le \varphi_x\left(x_0,t_0\right) \le p_l'.$$

Now since $p$ is a Riemann solution, either

$$0 \le p_r' \le p_l' \quad \text{or} \quad p_r' \le p_l' \le 0.$$

We first assume that $0 \le p_r' \le p_l'$, then by (**A**.**3**),

$$\begin{aligned} H\left(p_l',a_r\right) &\le H\left(\varphi_x\left(x_0,t_0\right),a_r\right) \le H\left(p_r',a_r\right) \\ &= H\left(p_l',a_l\right) \le H\left(\varphi_x\left(x_0,t_0\right),a_l\right) \le H\left(p_r',a_l\right). \end{aligned} \tag{3.25}$$

Also, since $u-\varphi$ has a maximum at $(0,t_0)$,

$$\frac{\varphi(0,t_0)-\varphi(0,t)}{t_0-t} \le \frac{u(0,t_0)-u(0,t)}{t_0-t},$$

for $t<t_0$ and $t_0-t$ sufficiently small. Now by definition

$$u(0,t_0) = u(t,0) - (t_0-t)\,H\left(p_l',a_l\right) = u(t,0) - (t_0-t)\,H\left(p_r',a_r\right),$$

and therefore

$$-H\left(p_r',a_r\right) = -H\left(p_l',a_l\right) \ge \frac{u(0,t_0)-u(0,t)}{t_0-t} \ge \frac{\varphi(0,t_0)-\varphi(0,t)}{t_0-t}.$$

Passing to the limit as $t\uparrow t_0$,

$$-H\left(p_r',a_r\right) = -H\left(p_l',a_l\right) \ge \varphi_t\left(0,t_0\right), \tag{3.26}$$

so

$$\varphi_t\left(0,t_0\right) + H\left(\varphi_x\left(0,t_0\right),a_r\right) \le -H\left(p_r',a_r\right) + H\left(p_r',a_r\right) = 0.$$

By (3.25) $\min\left\{H(\varphi_x(0,t_0),a_l), H(\varphi_x(0,t_0),a_r)\right\} = H(\varphi_x(0,t_0),a_r)$ and $u$ is a viscosity subsolution.

Now assume that $p_r' \le p_l' \le 0$, then we find

$$\begin{aligned} H\left(p_r',a_l\right) &\le H\left(\varphi_x\left(x_0,t_0\right),a_l\right) \le H\left(p_l',a_l\right) \\ &= H\left(p_r',a_r\right) \le H\left(\varphi_x\left(x_0,t_0\right),a_r\right) \le H\left(p_l',a_r\right). \end{aligned} \tag{3.27}$$

Then by (3.26)

$$\varphi_t\left(0,t_0\right) + H\left(\varphi_x\left(0,t_0\right),a_l\right) \le -H\left(p_l',a_l\right) + H\left(p_l',a_l\right) = 0,$$

and in this case, by (3.27), $\min\left\{H(\varphi_x(0,t_0),a_l), H(\varphi_x(0,t_0),a_r)\right\} = H(\varphi_x(0,t_0),a_l)$ so $u$ is a subsolution.

To verify that $u$ is also a supersolution, we can use the same arguments. □

In order to define front tracking for Hamilton–Jacobi equations, we now observe that the same formula, (3.23) holds if $H$ is replaced by the piecewise linear interpolant $H^\delta$. Let a front in $p^\delta$ start at a point $(x_0, t_0)$, and let its position be

$$x(t) = x_0 + \sigma\,(t - t_0)\,,$$

where $\sigma$ is the speed of the front. Now define $u^\delta(x(t), t)$ by

$$\begin{aligned} u^\delta(x(t),t) &= u^\delta\,(x_0,t_0) + (x - x_0)\,p^\delta(x(t)^\pm,t) - (t-t_0)\,H^\delta\left(p^\delta(x(t)^\pm, a^\delta\left(x(t)^\pm\right), g^\delta(t)\right) \\ &= u^\delta\,(x_0,t_0) + (t-t_0)\left(\sigma p^\delta(x(t)^\pm,t) - H^\delta\left(p^\delta(x(t)^\pm, a^\delta\left(x(t)^\pm\right), g^\delta(t)\right)\right). \end{aligned} \tag{3.28}$$

Between fronts we define $u^\delta$ by linear interpolation. By Lemma 3.3 it follows that $u^\delta$ is the unique viscosity solution to the initial value problem

$$u^\delta_t + H^\delta\left(u^\delta_x, a^\delta(x), g^\delta(t)\right) = 0, \quad u^\delta(x,0) = u^\delta_0(0) + \int_0^x p^\delta(x,0)\,dx, \tag{3.29}$$

where the piecewise linear interpolant $H^\delta$ is defined via the front tracking procedure, c.f., (3.7), and $u^\delta_0(0)$ is a constant. Furthermore

$$u^\delta(x,t) = u^\delta_0(0) - tH^\delta_0 + \int_0^x p^\delta(x,t)\,dx, \tag{3.30}$$

where $H^\delta_0 = H^\delta(p^\delta(0^\pm,t), a^\delta(0^\pm), g^\delta(t))$. Since $p^\delta$ converges to $p$ in $L^1$, and the traces $H(p^\delta, a^\delta, g^\delta)$ exists for almost all $t$, it follows that $u^\delta$ also converges in $L^\infty$, and its limit is

$$u(x,t) = u_0(0) - tH_0 + \int_0^x p(x,t)\,dx. \tag{3.31}$$

Now assume that $u - \varphi$ has a local maximum at $(x_0, t_0)$, then since $u^\delta$ converges in $L^\infty$ to $u$ $u^\delta - \varphi$ has a local maximum at

$$\left(x^\delta, t^\delta\right),$$

and $x^\delta \to x$ and $t^\delta \to t$ as $\delta \to 0$. Since $u^\delta$ is a viscosity (sub)solution

$$\begin{aligned} &\varphi_t\left(x^\delta, t^\delta\right) \\ &\quad + \min\left\{H^\delta\left(\varphi_x\left(x^{\delta,-}, t^{\delta,-}\right), a^\delta\left(x^{\delta,-}\right), g^\delta\left(t^{\delta,-}\right)\right), H^\delta\left(\varphi_x\left(x^{\delta,+}, t^{\delta,-}\right), a^\delta\left(x^{\delta,+}\right), g^\delta\left(t^{\delta,-}\right)\right)\right\} \le 0, \end{aligned}$$

Sending $\delta \downarrow 0$, we find that the limit $u$ is a viscosity subsolution. To prove that $u$ is also a supersolution, we can use analogous arguments. Thus $u$ is a viscosity solution of (1.1), with the initial data

$$u(x,0) = u_0(0) + \int_0^x p(x,0)\,dx.$$

Furthermore, for any fixed $y$ we have that

$$\begin{array}{ccccc} u^\delta(x,t) = & u^\delta(y,t) + & \displaystyle\int_y^x p^\delta(z,t)\,dz, & & \\ \downarrow & \downarrow & \downarrow & & \text{as } \delta \to 0, \\ u(x,t) = & u(y,t) + & \displaystyle\int_y^x p(z,t)\,dz, & & \end{array}$$

which means that $u_x = p$ almost everywhere. Thus we have proved existence of a viscosity solution, and that this solution can be constructed by front tracking, i.e.,

**Theorem 3.5.** *Assume that (**A**.1) – (**A**.8) all hold. For each $u_0 \in BUC(\mathbb{R})$ and all $T > 0$, there exists a viscosity solution $u \in BUC(\mathbb{R} \times [0,T])$ of* (1.1)*. Moreover*

$$p(x,t) = \frac{\partial u}{\partial x}(x,t),$$

*is a weak entropy solution of* (1.14)*, with the initial condition*

$$p(x,0) = \frac{\partial u_0}{\partial x}(x).$$

*The viscosity solution $u$ can be constructed by the front tracking scheme outlined above.*

## References

[1] F. E. Benth, K. H. Karlsen, and K. Reikvam, *Portifolio optimization in a Lévy market with intertemporal substitution and transaction costs*, Stoch. and stoc. reports, (2002).

[2] I. Capuzzo Dolcetta and B. Perthame, *On some analogy between different approaches to first order PDE's with nonsmooth coefficients*, Adv. Math. Sci. Appl., (1996).

[3] G. M. Coclite and N. H. Risebro, *Conservation laws with a time dependent discontinuous coefficients*, preprint, available at the url: `www.math.ntnu.no/conservation`, (2002).

[4] M. G. Crandall, H. Ishii, and P.-L. Lions, *User's guide to viscosity solutions to second order partial differential equations*, Trans. AMS, (1992).

[5] G. Dal Maso and H. Frankowska, *Value functions for Bolza problems with discontinuous Lagrangians and Hamilton-Jacobi inequalities*, ESAIM Control Optim. Calc. Var., 5 (2000), pp. 369–393 (electronic).

[6] G. Dal Maso and H. Frankowska, *Autonomous integral functionals with discontinuous nonconvex integrands: Lipschitz regularity of minimizers, DuBois-Reymond neccesary conditions, and Hamilton–Jacobi equations*, preprint, available at the url: `cvgmt.sns.it/papers/dalfra02`, (2002).

[7] T. Gimse and N. H. Risebro, *Riemann problems with a discontinuous flux function*, in Proc. 3rd Internat. Conf. Hyperbolic Problems, Uppsala, 1991, Studentlitteratur, pp. 488–502.

[8] ———, *Solution of the Cauchy problem for a conservation law with a discontinuous flux function*, SIAM J. Math. Anal., 23 (1992), pp. 635–648.

[9] H. Ishii and M. Ramaswami, *Uniqueness results for a class of Hamilton–Jacobi equations with singular coefficients*, Comm. Partial Differential Equations, (1995).

[10] K. H. Karlsen and N. H. Risebro, *A note on front tracking and the equivalence between viscosity solutions of Hamilton-Jacobi equations and entropy solutions of scalar conservation laws*, Nonlinear Analysis: Theory, Methods & Applications, 50 (2002), pp. 455–469.

[11] K. H. Karlsen, N. H. Risebro, and J. D. Towers, *$L^1$ stability for entropy solutions of nonlinear degenerate parabolic convection-diffusion equations with discontinuous coefficients.* Preprint, available at the url `www.math.ntnu.no/conservation`, 2003.

[12] R. A. Klausen and N. H. Risebro, *Stability of conservation laws with discontinuous coefficients*, J. Differential Equations, 157 (1999), pp. 41–60.

[13] C. Klingenberg and N. H. Risebro, *Convex conservation laws with discontinuous coefficients. Existence, uniqueness and asymptotic behavior*, Comm. Partial Differential Equations, 20 (1995), pp. 1959–1990.

[14] ———, *Stability of a resonant system of conservation laws modeling polymer flow with gravitation*, J. Differential Equations, 170 (2001), pp. 344–380.

[15] S. N. Kružkov, *First order quasi-linear equations in several independent variables*, Math. USSR Sbornik, 10 (1970), pp. 217–243.

[16] D. N. Ostrov, *Viscosity solutions and convergence of monotone schemes for synthetic aperture radar shape-from-shading equations with discontinuous intensities*, SJAM, 59 (1999), pp. 2060–2085.

[17] ———, *Extending viscosity solutions to Eikonal equations with discontinuous spatial dependence*, Nonlin. Anal., 42 (2000), pp. 709–736.

[18] N. Seguin and J. Vovelle, *Analysis and approximation of a scalar conservation law with a flux function with discontinuous coefficients.*, Preprint, (2002).

[19] T. Strömberg, *On viscosity solutions of irregular Hamilton–Jacobi equations*, tech. rep., 2003. Preprint.

[20] B. Temple, *Global solution of the Cauchy problem for a class of $2 \times 2$ non-strictly hyperbolic conservation laws*, Adv. in Appl. Math., 3 (1982), pp. 335–375.

[21] A. Tourin, *A comparison theorem for a piecewise Lipschitz continuous Hamiltonian and application to shape from shading problems*, Numer. Math., (1992).

[22] J. D. Towers, *Convergence of a difference scheme for conservation laws with a discontinuous flux*, SIAM J. Num. Anal., 38 (2000), pp. 681–698.

[23] ———, *A difference scheme for conservation laws with a discontinuous flux - the nonconvex case*, SIAM J. Num. Anal., 39 (2001), pp. 1197–1218.

(G. M. Coclite)
SISSA–ISAS
via Beirut 2–4
34014 Trieste
Italy
*E-mail address*: `coclite@sissa.it`

(N. H. Risebro)
Department of Mathematics
University of Oslo
P.O. Box 1053, Blindern
N–0316 Oslo, Norway
*E-mail address*: `nilshr@math.uio.no`